\documentclass[12pt]{article}
\usepackage[final]{epsfig}
\usepackage{graphics}
\usepackage{amsmath}
\usepackage{amsfonts}
\usepackage{latexsym}
\usepackage{amssymb}
\usepackage{graphicx}
\usepackage{url}
\usepackage{epstopdf}
\usepackage{wasysym}

\newtheorem{lemma}{Lemma}[section]

\newtheorem{remark}[lemma]{Remark}

\newtheorem{theorem}{Theorem}

\begin{document}
\newcommand{\eps}{{\varepsilon}}
\newcommand{\proofend}{$\Box$\bigskip}
\newcommand{\C}{{\mathbb C}}
\newcommand{\Q}{{\mathbb Q}}
\newcommand{\R}{{\mathbb R}}
\newcommand{\Z}{{\mathbb Z}}
\newcommand{\RP}{{\mathbb {RP}}}
\newcommand{\CP}{{\mathbb {CP}}}
\newcommand{\Tr}{\rm Tr}
\def\proof{\paragraph{Proof.}}

\title{Two variations on the periscope theorem}

\author{Serge Tabachnikov\footnote{
Department of Mathematics,
Penn State University,
University Park, PA 16802;
tabachni@math.psu.edu}
}

\date{}
\maketitle

\section{Introduction} \label{sec:intro}

Geometrical, or ray, optics, is a classical subject that remains an active area of research. 

One of the reasons for the contemporary interest in geometrical optics is a spectacular recent progress of freeform optical design. From a mathematical point of view, a freeform mirror is a smooth hypersurface in $\R^n$ that reflects oriented lines (rays of light) according to the familiar law ``the angle of incidence equals the angle of reflection". That is, the outgoing ray is the reflection of the incoming ray 
in the tangent hyperplane to the mirror at the impact point.

Geometrical optics is a source of important examples and open problems in symplectic geometry: the space of oriented lines ${\mathcal L}^{2n-2}$ in $\R^n$ carries a symplectic structure (symplectimorphic to the canonical symplectic structure of the cotangent bundle $T^*S^{n-1}$), and the reflection in a mirror is a local symplectomorphism of ${\mathcal L}^{2n-2}$.  

We refer to \cite{Bo,Ko} for panoramic views of geometrical optics and to \cite{Ma} for a modern approach to some classical results. The reader who is interested in the original treatment of this subject by W. R. Hamilton is referred to \cite{Ha}. 

Let $M^{n-1}$ be a (germ of a) cooriented hypersurface in $\R^n$. The normals to $M$ form a normal family of oriented lines, and $M$ is its front. The front is not unique: a normal family has a 1-parameter family of equidistant fronts. 
In terms of symplectic geometry, a normal family is  characterized as a Lagrangian submanifold in ${\mathcal L}^{2n-2}$.

Since the reflection in a mirror is a symplectic transformation, a normal family reflects to a normal family. This is the Malus-Dupin theorem. 

A converse statement is due to Levi-Civita \cite{LC}:  given two generic local normal families, consisting of the outgoing and the incoming rays, there is a one-parameter family of mirrors that reflect one family to the other. The mirrors are the loci of points for which  the sum of distances to the two wave fronts is constant (generalizing the gardener's, or string, construction of an ellipse). 

Thus, given two fronts, $M_1^{n-1}, M_2^{n-1} \subset \R^n$, there is a 1-parameter family of local diffeomorphisms defined by the optical reflections. It is an open problem to describe such diffeomorphisms of germs of hypersurfaces.

A 2-mirror system that reflects one normal family to another can be constructed by choosing the first mirror arbitrarily and the second one according to the Levi-Civita theorem. A mirror is locally a graph of a function of $n-1$ variables, so such a 2-mirror system depends on one functional parameter. 

More generally, an $n$-mirror system that reflects one normal family to another depends on $n-1$ functions of $n-1$ variables. Since a
local diffeomorphism of $(n-1)$-dimensional manifolds depends on $n-1$ functions of $n-1$ variables,  one  expects to need at least $n$ mirrors to realize a generic local diffeomorphism $M_1^{n-1} \to M_2^{n-1}$. 

To the best of our knowledge, the exact number of mirrors needed is not known. The case of $n=3$ was investigated in \cite{HC}; using the theory of exterior differential systems, it was shown that four mirrors sufficed.

\begin{figure}[hbtp]
\centering
\includegraphics[height=2.3in]{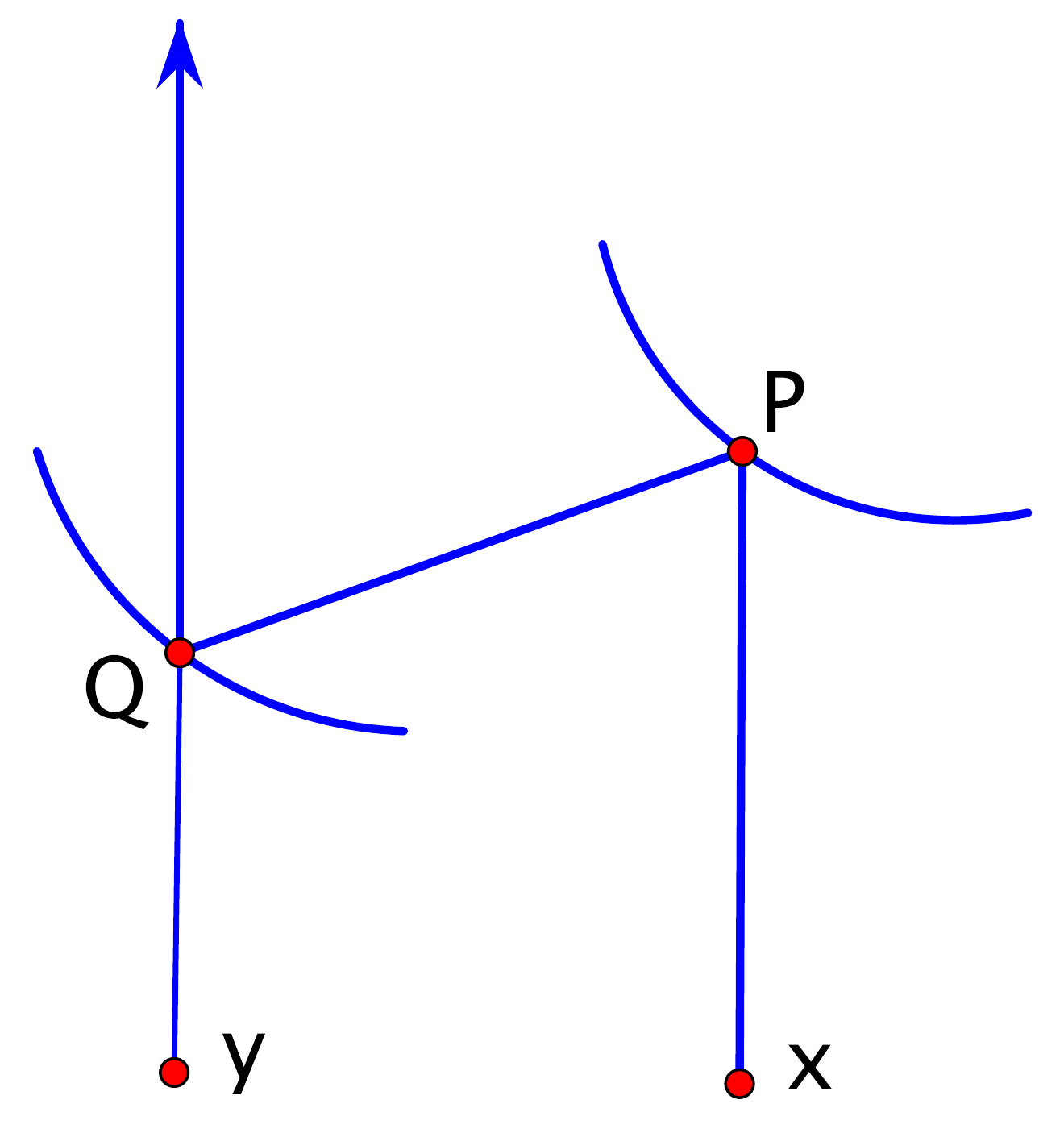} 
\caption{A periscope: a ray $xP$ reflects to the ray $yQ$.}
\label{figper}
\end{figure}

In the present paper, we are concerned by reflections in two mirrors, the case of one functional parameter. We present  two variations on the {\it periscope theorem} \cite{GO,O,PTT}. A periscope is a system of two mirrors that reflect the rays of light having a fixed (say, vertical) direction to the rays having the same direction. See Figure \ref{figper}. 

The periscope theorem states that the local diffeomorphism $x \mapsto y$ of the fronts, which are $(n-1)$-dimensional discs, is gradient: $y=x + \nabla f(x)$, where $f$ is a function of $n-1$ variables, depending on the mirrors. 

We consider two situations. A {\it spherical periscope} is a system of two mirrors that reflect the rays emanating from a fixed point $O$ to the rays coming back to $O$. A {\it reversed periscope} is a system of two mirrors that reflect the rays having a fixed direction to the rays having the opposite direction. 

In both cases, we describe the related local diffeomorphisms of the fronts (spherical, in the former, and flat, in the latter cases). For spherical periscopes, these are Theorems \ref{thm:main} and \ref{thm:express} in Section \ref{sec:sph}, and for reversed periscopes, this is Theorem \ref{thm:antiper} in Section \ref{sec:rev}.  

\begin{remark}
{\rm
Let us also mention a paper by R. Perline \cite{Per}, in which a somewhat related problem was studied: the optical reflection in thin films, that is, double mirror systems, in the limit as the two mirrors approach each other.
}
\end{remark}

\section{Spherical periscope} \label{sec:sph}

\begin{figure}[hbtp]
\centering
\includegraphics[height=2.3in]{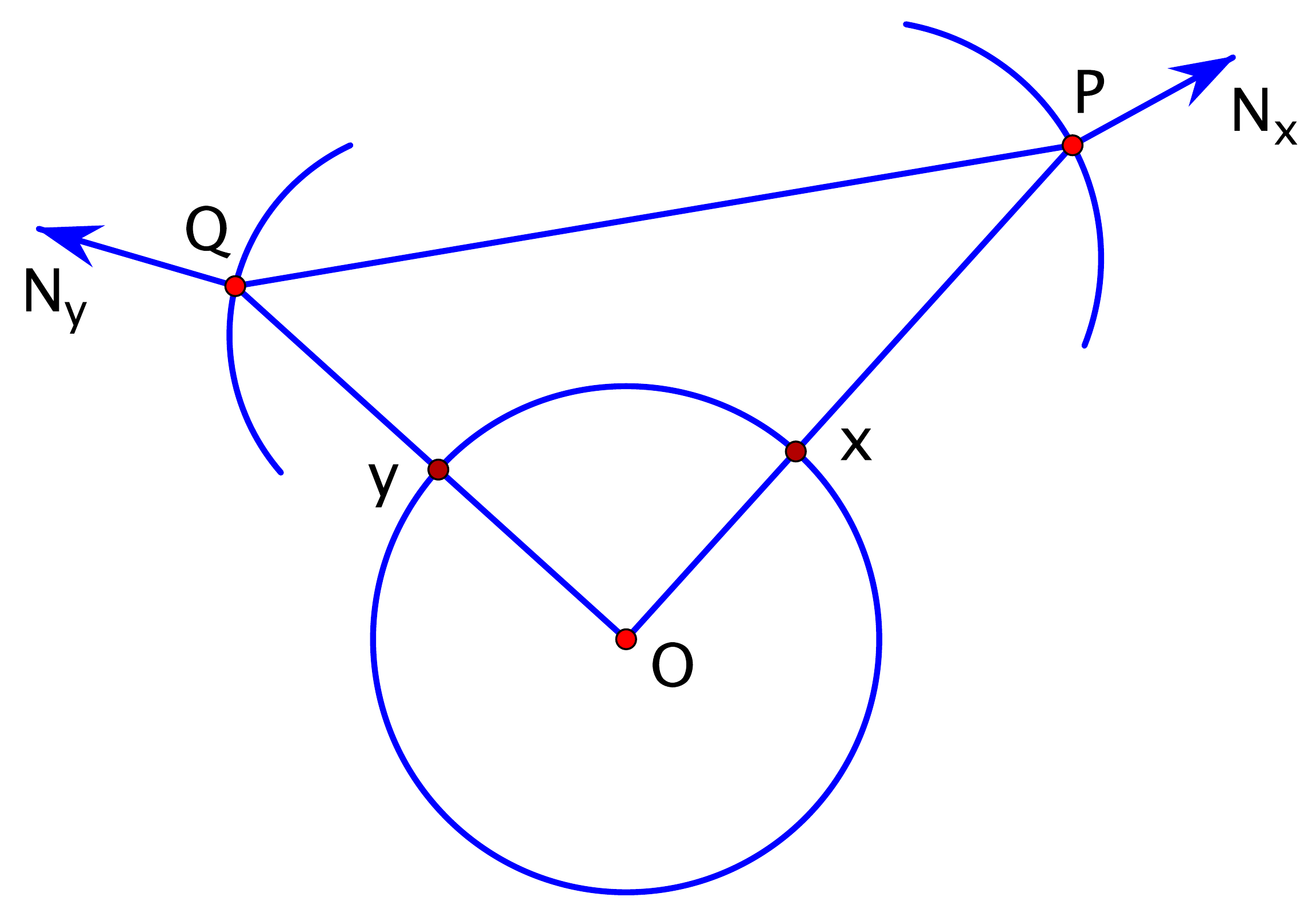} 
\caption{A spherical periscope.}
\label{mirrors}
\end{figure}

Consider the following situation: a ray of light $Ox$, emanating from point $O$, consecutively reflects in two mirrors and returns to point $O$ as the ray $yO$, see Figure \ref{mirrors}. Assume that the same holds for all rays, emanating from point $O$ and sufficiently close to the ray $Ox$, that is, for a neighborhood of point $x \in S^{n-1}$. %We refer to this system of two mirrors as a {\it spherical periscope}. 
One has a local diffeomorphism $T: x \mapsto y$ of the sphere $S^{n-1}$. 

Let $T$ be a local diffeomorphism of the sphere that does not take any point to the antipodal point. For every $x$ in the domain of $T$, connect $x$ with $T(x)$ by the shortest geodesic arc, and let $V_T (x)$ be the unit tangent vector to this arc at point $x$. This construction associates with $T$ a unit tangent vector field in its domain.

Call a vector field on a Riemannian manifold {\it projectively gradient} if it is proportional to a gradient vector field (with a functional factor).

A gradient vector $\nabla F$ field is orthogonal to the level hypersurface of the function $F$, so the normals to
a projectively gradient vector field for an integrable codimension 1 distribution. Replacing a vector field by its dual differential 1-form $\alpha$, the condition for being projectively gradient is $\alpha \wedge d\alpha =0$. In particular, the vector field, corresponding to a contact 1-form, is not projectively gradient. An example of such a field in $\R^3$ is
$$
y \frac{\partial}{\partial x} + \frac{\partial}{\partial z}.
$$

Our first result is as follows.

\begin{theorem} \label{thm:main}
Given a spherical periscope, the vector field $V_T$ is projectively gradient. 
\end{theorem}

\proof
Let $O$ be the origin and the sphere $S^{n-1} \subset \R^n$ be unit. 
Let us characterize the mirror hypersurfaces by their radial functions. Let $f: S^{n-1} \to \R$ be a (locally defined) smooth function, and let $P(x) = e^{f(x)} x$. Similarly, $Q(y) = e^{g(y)} y$ for another function $g: S^{n-1} \to \R$.

We claim that the vector $N_x := x - \nabla f(x)$ is normal to the first mirror at point $P(x)$. Indeed, let $v \in T_x S^{n-1}$ be a test tangent vector. Then
$$
\lim_{\eps \to 0}\frac{d \left[e^{f(x+\eps v)} (x+\eps v)\right]} {d \eps} = e^{f(x)} [v + (v \cdot \nabla f(x)) x],
$$
and
$$
[v + (v \cdot \nabla f(x)) x] \cdot [x - \nabla f(x)] = 0,
$$
since $x \cdot v =0, x \cdot x =1$, and $x \cdot \nabla f(x) =0$.

Similarly, $N_y := y - \nabla g(y)$ is normal to the second mirror at point $Q$. 

Next we claim that the vectors $x,y, PQ, N_x$, and $N_y$ lie in the same 2-plane. Indeed, let $\pi$ be the plane $OPQ$. Then the first three vectors obviously lie in $\pi$. By the definition of mirror reflection, the normal $N_x$ is coplanar with the incoming and outgoing rays, hence $N_x \in \pi$, and likewise for $N_y$.

It follows that the projection of $N_x$ to the tangent hyperplane $T_x S^{n-1}$ is $-\nabla f(x)$, it lies in the plane $\pi$ and it is tangent to the geodesic arc $xy$ at point $x$.  
\proofend

Next we calculate, in terms of the function $f(x)$, the other function, $g(y)$, and the spherical distance $d(x,y)$ between points $x$ and $y$.

Let $\alpha$ be the (acute) angle between the vectors $x$ and $N_x$, and $\beta$ that between $y$ and $N_y$. We abbreviate $f(x)$ and $g(y)$ to  $f$ and $g$.

\begin{lemma} \label{lm:angles}
One has:
\begin{equation} \label{eq:angles}
\tan \alpha = |\nabla f|,\ \tan \beta = |\nabla g|.
\end{equation}
\end{lemma}

\proof
One has:
$$
\cos\alpha = \frac{x \cdot N_x}{|N_x|} = \frac{1}{\sqrt{1+ |\nabla f|^2}}.
$$
This implies the formula for $\tan \alpha$, and similarly for $\tan \beta$.
\proofend

\begin{lemma} \label{lm:rel1}
One has
\begin{equation} \label{eq:rel1}
\frac{e^f |\nabla f|}{1+ |\nabla f|^2} = \frac{e^g |\nabla g|}{1+ |\nabla g|^2}.
\end{equation}
\end{lemma}

\proof
The sine rule in triangle $OPQ$ implies that $e^f \sin (2\alpha) = e^g \sin (2\beta)$. Expressing the sine of the double angle via the tangent of the angle and using Lemma \ref{lm:angles} yields the result.
\proofend

Let $S(x,y)$ be the common value of (\ref{eq:rel1}).

\begin{lemma} \label{lm:rel2}
One has
$$
S = C \frac{|\nabla f| |\nabla g|}{|\nabla f| + |\nabla g|},
$$
where $C$ is a constant.
\end{lemma}

\proof
The perimeter of triangle $OPQ$ is constant (that is, does not depend on the point $x$). This is due to the fact well known in geometrical optics:  the optical path length  between two wave fronts is the same for all rays\footnote{The length function is constant on its critical manifold.}; in our situation, both wave fronts are subsets of the unit sphere $S^{n-1}$. Denote this perimeter by $2C$, that is, $|PQ| = 2C - e^f - e^g$. 

The cosine rule in triangle $OPQ$ implies that
$$
e^{2f} + e^{2g} - 2e^{f+g} \cos (\pi-2\alpha-2\beta) = (2C - e^f - e^g)^2, 
$$
or
\begin{equation} \label{quad1}
e^{f+g} [1 - \cos (2\alpha+2\beta)] - 2C (e^f + e^g) + 2C^2 =0.
\end{equation}

We express $\cos(2\alpha), \sin(2\alpha), \cos(2\beta), \sin(2\beta)$ via $\tan\alpha$ and $\tan \beta$, and then, using formulas (\ref{eq:angles}), via $|\nabla f|$ and $|\nabla g|$. This results in
\begin{equation} \label{eq:cos}
1 - \cos (2\alpha+2\beta) = \frac{2 (|\nabla f| + |\nabla g|)^2}{(1+ |\nabla f|^2) (1+ |\nabla g|^2)}.
\end{equation}

Using formulas (\ref{eq:rel1}), we  express $e^f$ and $e^g$ in terms of $|\nabla f|, |\nabla g|$, and $S$. Substituting this in equation (\ref{quad1}) yields a quadratic equation in $S$:
$$
S^2 (|\nabla f| + |\nabla g|)^2 - CS (|\nabla f| + |\nabla g|)(1+ |\nabla f| |\nabla g|) + C^2 |\nabla f| |\nabla g| = 0.
$$

This equation has two solutions:
$$
S_1=\frac{C}{|\nabla f| + |\nabla g|}\ \ {\rm and}\ \ S_2=\frac{C |\nabla f| |\nabla g|}{|\nabla f| + |\nabla g|}.
$$
To see that we should select the second one, we argue as follows. 

If $f$ is constant, then the first mirror is a sphere centered at the origin. Then the second mirror is also spherical, and $\nabla f = \nabla g =0$. Then expression (\ref{eq:rel1}) vanishes, but $S_1$ has infinite value. Hence $S_1$ is an extraneous root.
\proofend

Now we can express the function $g$ and the spherical distance $d(x,y)$  in terms of $f$.

\begin{theorem} \label{thm:express}
One has
$$
e^g = \frac{e^{2f} - 2C e^f + C^2(1 + |\nabla f|^2)}{C(1+ |\nabla f|^2) - e^f},
$$
and
$$
d(x,y) = \pi - 2 \arcsin \left[ \frac{C |\nabla f|}{\sqrt{e^{2f} - 2C e^f + C^2(1 + |\nabla f|^2)}}   \right].
$$
\end{theorem}

\proof
First, we determine $|\nabla g|$ from the equation
$$
\frac{e^f |\nabla f|}{1+ |\nabla f|^2} = C \frac{|\nabla f| |\nabla g|}{|\nabla f| + |\nabla g|};
$$
we get
\begin{equation} \label{eq:nabla}
|\nabla g| = \frac{e^f |\nabla f|}{C(1+ |\nabla f|^2) - e^f}.
\end{equation}
Therefore 
\begin{equation} \label{eq:nabla2}
1+ |\nabla g|^2 = \frac{[1+ |\nabla f|^2] [e^{2f} - 2C e^f + C^2(1 + |\nabla f|^2)]}{[C(1+ |\nabla f|^2) - e^f]^2}.
\end{equation}
Now, using equation (\ref{eq:rel1}), we find 
$$
e^g = \frac{e^{2f} - 2C e^f + C^2(1 + |\nabla f|^2)}{C(1+ |\nabla f|^2) - e^f},
$$
as claimed.

Next, we use equation (\ref{eq:cos}):
$$
\sin(\alpha+\beta) = \sqrt {\frac{1 - \cos (2\alpha+2\beta)}{2}} = \frac{(|\nabla f| + |\nabla g|)}{(\sqrt{1+ |\nabla f|^2}) (\sqrt{1+ |\nabla g|^2)}}.
$$
Substitute $|\nabla g|$ from (\ref{eq:nabla})  and $1+ |\nabla g|^2$ from (\ref{eq:nabla2})
to obtain
$$
\sin(\alpha+\beta) = \frac{C |\nabla f|}{\sqrt{e^{2f} - 2C e^f + C^2(1 + |\nabla f|^2)}}. 
$$
Since $d(x,y) = \pi - 2(\alpha + \beta)$, this implies the result.
\proofend

\section{Reversed periscope} \label{sec:rev}

\begin{figure}[hbtp]
\centering
\includegraphics[height=2in]{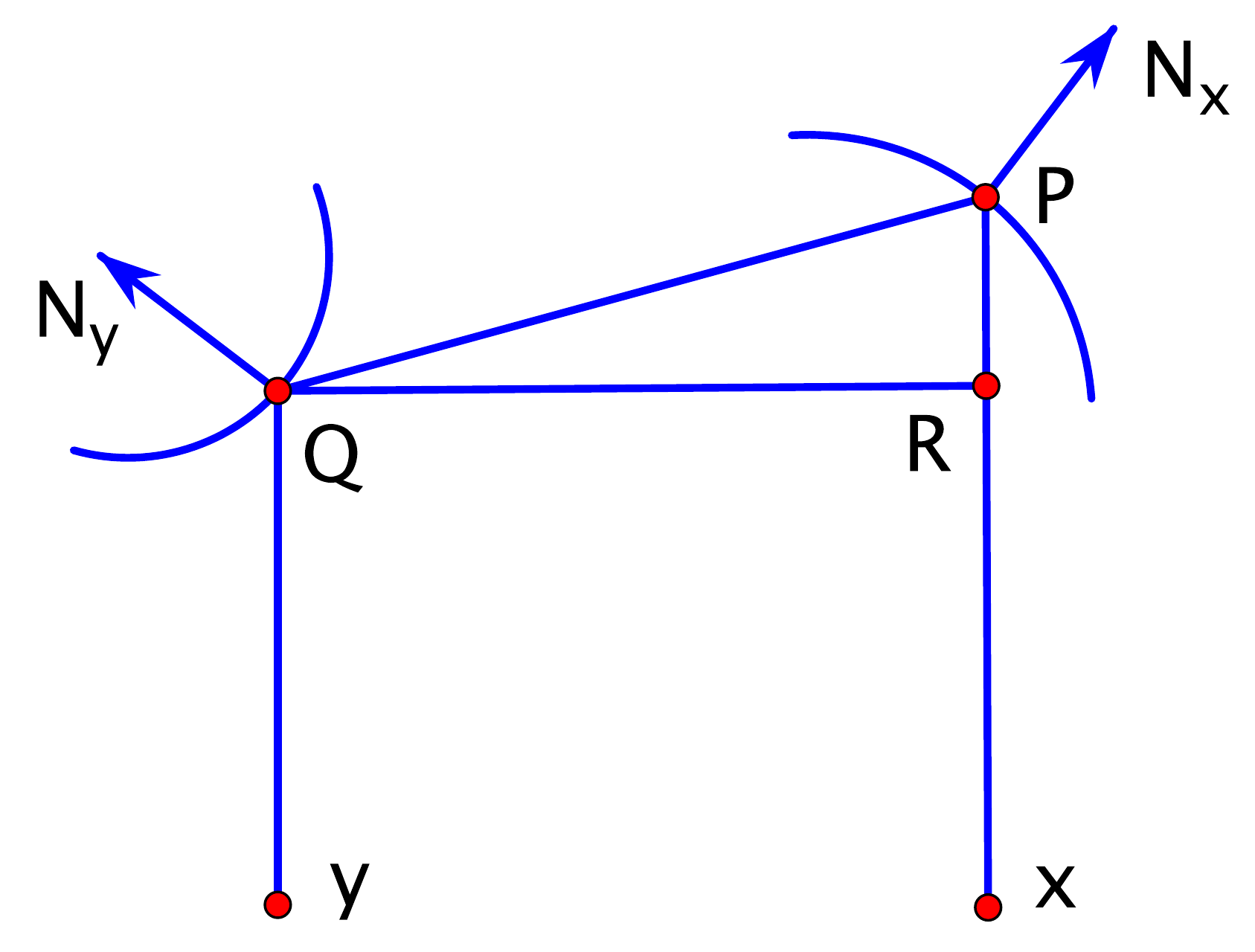} 
\caption{A reversed periscope: the ray $xP$ reflects to the ray $Qy$.}
\label{antiper}
\end{figure}

We now consider the situation similar to the previous one, but with the point $O$ located at infinity. More precisely,  consider $\R^n$ as the product $\R^{n-1} \times \R$ and think of the last coordinate axis as vertical. We have two mirrors and every upward vertical ray of light (in a certain neighborhood of one such ray) reflects to a downward ray.   See Figure \ref{antiper}.

Vertical rays are parameterized by points of the horizontal hyperplane $\R^{n-1}$, and we have a local diffeomorphism $T: x \mapsto y$. Let $U(x)$ be the vector $xy$.

Without loss of generality, we assume that point $P$ is not lower than point $Q$: otherwise, we reverse the directions of the rays and interchange $x$ and $y$. The mirrors are graphs of (locally defined) functions $f(x): \R^{n-1} \to \R$ and $g(y): \R^{n-1} \to \R$. 

The next result is an analog of Theorem \ref{thm:express}: it expresses the function $g$ and the vector $T$ in terms of the function $f$.

\begin{theorem} \label{thm:antiper}
One has: 
\begin{equation} \label{eq:gbyf}
g = \frac{f - C(1-|\nabla f|^2)}{|\nabla f|^2},
\end{equation}
and
$$
T (x) = x + \frac{2(C-f(x))}{|\nabla f(x)|^2} \nabla f(x),
$$
where $C$ is a positive constant.
\end{theorem}

In particular, the vector field $U(x)$ is projectively gradient (an analog of Theorem \ref{thm:main}). 

\proof
The arguments are similar to the ones in Section \ref{sec:sph}. 

The normals to the mirrors at points $P$ and $Q$ are given by the formulas
$$
N_x = (-\nabla f, 1),\ N_y = (-\nabla g, 1).
$$
Let angle $RPQ$ be $2\alpha$. Then 
$$
\cos\alpha = \frac{N_x \cdot (0,1)}{|N_x|} = \frac{1}{\sqrt{1+ |\nabla f|^2}}.
$$
Hence $\tan\alpha = |\nabla f|$. Likewise,  
$$
\cos\left(\frac{\pi}{2}-\alpha\right) = \sin\alpha = \frac{1}{\sqrt{1+ |\nabla g|^2}}.
$$
Hence $\cot\alpha = |\nabla g|$, and therefore $|\nabla f| |\nabla g| =1$.

From the right triangle $PQR$, we have $|U|=(f-g) \tan (2\alpha)$, and using the formula for tangent of a double angle, we get
\begin{equation} \label{eq:mag}
|U|= \frac{2 (f-g) |\nabla f|}{1- |\nabla f|^2}.
\end{equation}
Since point $P$ is higher than point $Q$, we have $|\nabla f| < 1$.

Next, we use the fact that the optical path length is constant: 
$$
f+g + |PQ| = 2C,
$$
 and hence
$$
f+g+\frac{f-g}{\cos 2\alpha} =2C,
$$
or
$$
f(1+\cos 2\alpha) - g(1-\cos 2\alpha) = 2C \cos 2\alpha.
$$
Expressing everything in terms of $\tan \alpha = |\nabla f|$, yields formula (\ref{eq:gbyf}).

Finally, the vectors $U$ and $\nabla f$ lie in the plane $PQR$, hence $U$ is proportional to $\nabla f$, and due to our assumptions, with a positive coefficient. It follows from (\ref{eq:mag}) and (\ref{eq:gbyf}) that 
$$
U = \frac{2(C-f(x))}{|\nabla f(x)|^2} \nabla f(x),
$$
as claimed. 
\proofend

\bigskip
{\bf Acknowledgements}. Many thanks to P. Albers and A. Plakhov for valuable discussions. This work was supported by NSF grant DMS-1510055.

\end{document}